\documentclass[11pt]{article} 
\usepackage{amsmath,amsfonts,amssymb,color}
\usepackage{graphicx}
\usepackage[mathscr]{euscript}
 \usepackage{setspace}
\newcommand{\eps}{\epsilon}
\newcommand{\al}{\alpha}
\newcommand{\ov}{\overline}
\newcommand{\dis}{\displaystyle}
\newcommand{\pr}{\prime}
\newcommand{\R}{\mathbb R}

\newcommand{\dotp}{\mbox{\parbox[c]{.05in}{\huge{.}}} \, }

\begin{document}   
\title{Periodic Orbits for a Discontinuous Vector Field Arising from a Conceptual Model of Glacial Cycles}
\author{James Walsh\\Department of Mathematics, Oberlin College
\and Esther Widiasih\\Department of Mathematics, University of Hawaii--West Oahu
\and Jonathan Hahn\\ School of Mathematics, University of Minnesota
\and Richard McGehee\\ School of Mathematics, University of Minnesota}
\date{July  20, 2015}

\maketitle

\vspace{.25in}
\noindent{\bf Abstract.} \ Conceptual climate models provide an approach to understanding climate processes through a mathematical analysis of an approximation to reality.  Recently, these models have also provided interesting examples of nonsmooth dynamical systems.  Here we discuss a conceptual model of glacial cycles consisting of a system of three ordinary differential equations defining a discontinuous vector field.  We show that this system has a large periodic orbit crossing the discontinuity boundary.  This orbit can be interpreted as an intrinsic cycling of the Earth's climate giving rise to alternating glaciations and deglaciations.

\bigskip



\section{Introduction} 

In recent years an extensive theory of discontinuous differential equations has been built, based largely on the work of Filippov \cite{fil}.  Although the general theory is far from complete \cite{guck}, it is sufficiently developed to allow it to be applied to systems such as the model of glacial cycles described here.

Conceptual models of the Earth's ice sheets began with the work of Budyko \cite{bud} and Sellers \cite{sell}.  Using a simple equation incorporating latitudinal variation of the incoming solar radiation and the ice-albedo feedback effect, Budyko was able to show the existence of a stable small ice cap.  Widiasih \cite{esther} introduced an equation for the motion of the edge of the ice sheet (called the ice line), enabling the use of dynamical systems theory to examine the ice line dynamics.  Subsequently, she and McGehee \cite{dickesther}  introduced an approximation to the Budyko-Widiasih model, reducing the model to a system of two ordinary differential equations.

Many aspects of the Earth's glacial cycles remain a mystery.  For roughly the last million years, the Earth's climate has experienced extensive cold periods with large ice sheets covering much of North America and Eurasia.  These glacial periods have been followed by relatively rapid transitions to warmer temperatures as the ice sheets retreat.  These so-called interglacial periods, which are shorter in duration, are typically followed by a slow descent into the subsequent glacial period, leading to the asymmetric sawtooth pattern evident in the paleoclimate data \cite{petit}.  These glacial/interglacial cycles have occurred with a periodicity of roughly 100,000 years.

The explanation for the 100,000 year cycles usually appeals in some way to the work of Milankovitch \cite{milank}, who examined the impact of variations in the Earth's orbital parameters on the incoming solar radiation (called Milankovitch cycles).  It turns out that the eccentricity of the Earth's orbit changes with a period of about 100,000 years, corresponding with the glacial cycles of the last million years.  However, the paleoclimate record shows that, for the period between five million and one million years ago, the glacial cycles occurred with a period of about 41,000 years, corresponding instead to changes in the Earth's obliquity (axial tilt).  Since the Milankovitch cycles have not changed significantly over the last five million years, it is clear that the Earth's climate cycles cannot be explained as a simple response to the Milankovitch forcing.

Many models have been proposed to explain the response of the Earth's climate to the Milankovitch cycles.  One of the earliest was in a paper of Maasch and Saltzman \cite{maasch}, where they proposed that, for the last million years, the Earth's climate system has had an intrinsic oscillation of 100,000 years, independent of the Milankovitch cycles.  Dynamically speaking, the autonomous system has a stable periodic orbit with a period of about 100,000 years.  When the Milankovitch cycles are introduced, the system becomes nonautonomous and falls into synchronization with the eccentricity cycles.  The change that occurred one million years ago can be explained as slow changes in the parameters of the Earth's climate system creating a Hopf bifurcation in the autonomous system, changing the response to the forcing.

Many other conceptual models have been introduced to study nonlinear climate feedbacks as they relate to the glacial cycles. These models often include a threshold that, once exceeded, triggers a switch from, say, a glacial to an interglacial period. In previous studies the threshold variable has been ice volume \cite{huy2007}, \cite{huy2011}; a combination of Milankovitch forcing and ice volume \cite{paill1998}, \cite{parr2003};  a function of ice volume, Milankovitch forcing and area of the Antarctic ice sheet \cite{paillpar2004}; and a switch connected to the rapid formation of sea ice \cite{tzip2003}. In related work, a relaxation model is subjected to astronomical forcing to produce oscillations in global ice volume in  \cite{ashwin}. In many of these works, model parameters are adjusted in attempts to generate output that qualitatively mirrors the paleoclimate data.

Starting with the McGehee-Widiasih version of Budyko's equation  \cite{dickesther} and inspired by the work of K\"all\'en et al \cite{kallen}, we consider Widiasih's ice line variable to be the albedo line, where the snow cover remains throughout the year, and introduce a new variable, corresponding to the extent of the ice sheet.  The thinking here is that, during the glacial advance, the snow is accumulating (forming the glacier) at a lower latitude than the ice sheet edge, and hence that the albedo line is at a lower latitude than the ice sheet edge.  On the other hand, the work of Wright and Stefanova \cite{wright}  indicates that, during the glacial retreat, forests growing on top of the retreating ice sheets create an albedo line at a higher latitude than the edge of the ice sheet.  These phenomena give rise to what we think of as two distinct climate regimes, one occurring during the glacial advance and the other occurring during the glacial retreat.  The boundary between these two regimes corresponds to a discontinuity in the vector field, allowing us to view the model as a Filippov system. 

Our model is a system of three ordinary differential equations with a discontinuity boundary consisting of a plane in the three-dimensional state space.  Using techniques reminiscent of singular perturbation theory, we show the existence of a large periodic orbit crossing the discontinuity boundary.

In terms of the theory of glacial cycles, this work represents only the beginning, analogous to the establishment of a periodic orbit in the Maasch and Satzman model.  Missing is an analysis of the response to the Milankovitch forcing.  The forcing will be easy to introduce into the model, since it is based on the Budyko model, where the Milankovitch cycles can be easily incorporated \cite{dickclar}.

In terms of the mathematical theory, there may be some work to be done in studying the response of a Filippov system to time variation of the parameters.  In addition, our analysis points to the need to develop a theory of singular perturbations for Filippov systems.

The paper is organized as follows. In the next section we introduce Budyko's model and its coupling with a dynamic ice line. In Section 3 we review the quadratic approximation introduced in \cite{dickesther}, which we then augment by including a snow line independent of the edge of the ice sheet, as presented in Section 4. This in turn leads to a discontinuous three-dimensional vector field, the dynamics of which we analyze in Section 5. We point to several problems of mathematical interest generated by the study of the model and we summarize this work in the concluding section.

\section{A coupled temperature-ice line model} 

Budyko's energy balance model concerns the average annual temperature in latitudinal zones in a world assumed to be symmetric about the equator. While Budyko focused solely on equilibrium   temperature distributions, we consider the time-dependent equation   \cite{heldsuarez}
\begin{equation}\label{budyko}
R\frac{\partial T(y,t)}{\partial t}=Qs(y)(1-\alpha (y))-(A+BT)-C(T-\ov{T}).
\end{equation}
In equation \eqref{budyko}, $y$ is the sine of the latitude, chosen  for convenience (the area of an infinitesimal latitudinal strip at $y$ is proportional to $dy$).
Note that due to symmetry considerations, $y\in[0,1]$, with $y=0$ the equator and $y=1$ the North Pole.  

The function $T=T(y,t)$  ($^\circ$C)  is the annual mean surface temperature on the  circle of latitude at $y$, and $R$ is the heat capacity of the Earth's surface
(with units  J/(m$^{2} \,  ^\circ$C)). The left-hand side of \eqref{budyko} represents the change in energy stored in the Earth's surface at  $y$; units on each side of \eqref{budyko} are Watts per meter squared (W/m$^{2}$).

$Q$ denotes the mean annual incoming solar radiation (or {\em insolation}), a parameter depending on the eccentricity of Earth's orbit \cite{dickclar}. The function $s(y)$, which depends upon the obliquity of Earth's orbit \cite{dickclar}, accounts for the distribution of insolation across latitude, and satisfies
\begin{equation}\label{sofy}\notag
\int^1_0 s(y) dy =1.
\end{equation} 
The function $\al(y)$ represents the planetary {\em albedo}, a measure of the extent to which insolation is reflected back into space. Thus the $Qs(y)(1-\alpha (y))$-term represents the energy from the sun absorbed at the surface at ``latitude" $y$.

The $(A+BT)$-term models the {\em outgoing longwave radiation} (OLR) emitted by the Earth, while the transport of heat energy across latitudes
is modeled by the $C(T-\ov{T})$-term, in which 
\begin{equation}\label{tbar}\notag
\ov{T}=\int^1_0 T(y,t) dy
\end{equation}
 is the global annual mean surface temperature. \ The constants $A, B$ and $C$ are positive and empirical (for a more detailed introduction to Budyko's EBM see \cite{dickesther}, \cite{jimchris}). 

In  \cite{dickesther} the equilibrium temperature profiles of \eqref{budyko} are shown to be given by
\begin{equation}\label{Tstarbud}
T^*(y)=\frac{1}{B+C}\left(Qs(y)(1-\al(y))-A+\frac{C}{B}(Q(1-\ov{\al})-A)\right),
\end{equation}
where 
\begin{equation}\label{albar}
\ov{\al}=\int^1_0 \al(y)s(y) dy.
\end{equation}
Letting $p_0(y)=1$ and $p_2(y)=\frac{1}{2}(3y^2-1)$ denote the first two even Legendre polynomials, we make  use of the expression 
\begin{equation}\label{approxsofy}
s(y)=s_0p_0(y)+s_2p_2(y),  \ s_0=1, \ s_2=-0.482, 
\end{equation}
which is uniformly within 2\% of the actual values of $s(y)$  \cite{north}, in all that follows.

In Budyko's model one assumes the planet has an ice cap, with ice at all latitudes above a certain latitude $y=\eta$, and no ice  south of $y=\eta$. The edge of the ice sheet $\eta$ is called the {\em ice line}. The albedo function is then given by 
\begin{equation}\label{budalb}
\al_\eta(y)=
 \begin{cases}
\al_1, &\text{if  \ } y<\eta \\
\al_2, &\text{if  \ } y>\eta,\\
\end{cases} \quad \quad
\al_1<\al_2,
\end{equation}
where $\al_1$ and $\al_2$ represent the albedos of surface with no ice cover and that having an ice cover, respectively. With this choice of albedo function, and using expression \eqref{approxsofy}, equilibrium temperature profiles \eqref{Tstarbud} are even, piecewise quadratic functions having  a discontinuity at $\eta$. We also note $\eta$ serves to parametrize \eqref{albar}, and hence equilibrium functions \eqref{Tstarbud}, as well. In particular, there are infinitely many equilibrium temperature functions, one for each value of $\eta$. We write $T^*_\eta(y)$ when we wish to emphasize this parametrization. 

Defining $T^*_\eta(\eta)=\frac{1}{2}(\lim_{y\rightarrow \eta^-} T^*(y)+\lim_{y\rightarrow \eta^+} T^*(y))$, one finds 
the temperature at equilibrium at the ice line is given by 
\begin{equation}\label{Tstarateta}
T^*_\eta(\eta)=\frac{1}{B+C}\left(Qs(\eta)(1-\al_0)-A+ \frac{C}{B}(Q(1-\ov{\al}(\eta)-A)\right),
\end{equation}
with $\al_0=\frac{1}{2}(\al_1+\al_2)$. Budyko was interested in the existence of $\eta$-values for which $T^*_\eta(\eta)=T_c$, where $T_c$ is a {\em critical temperature} above which ice melts and below which ice forms. Specifically, he investigated  the relationship between the number of such $\eta$-values and the parameter $Q$.

As mentioned in the introductory section, Budyko was motivated by the study of the positive ice-albedo feedback: If temperatures were to decrease, an existing ice sheet would grow, increasing the albedo and further lowering temperatures, leading to an ever larger ice sheet. A warming climate would lead to a smaller ice sheet and reduced albedo, thereby raising temperatures and  further reducing the size of the ice sheet. 

Budyko found that for a range of $Q$-values there exist $\eta_1<\eta_2$ with $T^*_{\eta_i}(\eta_i)=T_c=-10^\circ$C, $i=1, 2,$ with $\eta_1$ in lower latitudes and $\eta_2$ nearer the North Pole. If $Q$  decreased sufficiently, however, there were no such $\eta$-values, with the implication being  that  the ice line would be positioned at the equator  for these $Q$-values (a so-called {\em Snowball Earth} event). This was viewed as a consequence of ice-albedo feedback. Notably, previous analyses of Budyko's model lacked any treatment of the stability of the ``preferred" equilibrium temperature distributions having  $T^*_\eta(\eta)=T_c$  (or, indeed, of $T^*_\eta(y)$ for any $\eta$) from a dynamical systems perspective, that is, in the infinite-dimensional setting intrinsic to \eqref{budyko}. 
 
Additionally,  Budyko's model lacks any mechanism by which the ice line $\eta$ is allowed to respond to changes in temperature.  This limitation was remedied by E. Widiasih in \cite{esther} through the addition of an ODE for the evolution of $\eta$, leading to the integro differential system
\begin{subequations}\label{budwid} 
\begin{eqnarray}
R\frac{\partial T}{\partial t} & =   & Qs(y)(1-\al(y,\eta))-(A+BT)-C(T-\overline{T})  \label{budwidA}  \\ 
\frac{d\eta}{dt} & =   &  \rho (T(\eta,t)-T_c).\label{budwidB}
\end{eqnarray}
\end{subequations}
Here, $\rho>0$ is a parameter governing the relaxation time of the ice sheet. The temperature distribution $T(y,t)$ evolves according to Budyko's equation \eqref{budwidA}, while the dynamics of $\eta$ are determined by the temperature at the ice line, relative to the critical temperature. The ice sheet retreats toward the pole if $T(\eta,t)>T_c$, and moves equatorward if $T(\eta,t)<T_c$. 

Working in the infinite-dimensional setting and with parameters as in Budyko's work, Widiasih proved the existence of equilibrium temperature--ice line pairs $(T_{\eta_i}(y),\eta_i), i=1, 2$, for $\rho $ sufficiently small. The ice line $y=\eta_2$ corresponds to a stable, small  ice cap, while  ice line $y=\eta_1$ corresponds to a large, unstable ice sheet (in particular, there are no oscillations of the ice sheet in the model). This work provided Budyko's pioneering and influential model
with a modern  dynamical systems perspective.

\vspace*{0.2in}
\section{The approximation of McGehee and Widiasih}

Recall  equilibrium solutions of Budyko's equation \eqref{budyko} are even and piecewise quadratic, with a discontinuity at $\eta$, when using albedo function \eqref{budalb} and expression \eqref{approxsofy}. This motivated the introduction of a quadratic approximation to system \eqref{budwid} in \cite{dickesther}. We outline this approach here, referring  the reader to \cite{dickesther} for all computational details.

Let $X$ denote the space of even, piecewise quadratic functions having a discontinuity at $\eta$. The four-dimensional linear space $X$ can be parametrized by the new variables $w_0, z_0, w_2$ and $z_2$ by letting
\begin{equation}\label{quadT}
T(y)=
\begin{cases}
w_0+\frac{1}{2}z_0+(w_2+\frac{1}{2}z_2)p_2(y), & y<\eta\\
w_0-\frac{1}{2}z_0+(w_2-\frac{1}{2}z_2)p_2(y), & y>\eta\\
w_0+w_2p_2(\eta), & y=\eta.
\end{cases}
\end{equation}
The choice of $T(\eta)$ in \eqref{quadT} coincides with $T(\eta)=\frac{1}{2}(\lim_{y\rightarrow \eta^-} T(y)+\lim_{y\rightarrow \eta^+} T(y))$, to maintain consistency with \eqref{Tstarateta}. One can show  that, in these variables,  
\begin{equation}\label{tbarA}
\ov{T}=\int^1_0 T(y) dy=w_0+z_0(\eta-\textstyle{\frac{1}{2}})+z_2P_2(\eta), \ \mbox{where } \ P_2(\eta)=\dis \int^\eta_0 p_2(y) dy.
\end{equation}

Plugging \eqref{quadT} and \eqref{approxsofy} into \eqref{budyko} and equating the coefficients of $p_0(y)$ (that is, the constant terms) and $p_2(y)$, respectively, yields the four equations 
\begin{subequations}\label{dick4}
\begin{align}
R(\dot{w}_0+\textstyle{\frac{1}{2}}\dot{z}_0)&=Q(1-\al_1)-A-(B+C)(w_0+\textstyle{\frac{1}{2}}z_0)+C\ov{T}\label{dick4a} \\
R(\dot{w}_0-\textstyle{\frac{1}{2}}\dot{z}_0)&=Q(1-\al_2)-A-(B+C)(w_0-\textstyle{\frac{1}{2}}z_0)+C\ov{T}\label{dick4b} \\
R(\dot{w}_2+\textstyle{\frac{1}{2}}\dot{z}_2)&=Qs_2(1-\al_1)-(B+C)(w_2+\textstyle{\frac{1}{2}}z_2) \label{dick4c} \\
R(\dot{w}_2-\textstyle{\frac{1}{2}}\dot{z}_2)&=Qs_2(1-\al_2)-(B+C)(w_2-\textstyle{\frac{1}{2}}z_2) \label{dick4d}
\end{align}
\end{subequations}
(two equations  each for $y<\eta$ and $y>\eta$). 
Adding and subtracting equations \eqref{dick4a}-\eqref{dick4b} and \eqref{dick4c}-\eqref{dick4d} and substituting in expression \eqref{tbarA} for $\ov{T}$, one arrives at
\begin{align}  \label{4ODE} 
R\dot{w}_0&= Q(1-\al_0)-A-Bw_0+C((\eta-\textstyle{\frac{1}{2}})z_0+z_2P_2(\eta))\\  \notag
R\dot{z}_0&=Q(\al_2-\al_1)-(B+C)z_0\\   \notag
R\dot{w}_2&=Qs_2(1-\al_0)-(B+C)w_2 \\  \notag
R\dot{z}_2&=Qs_2(\al_2-\al_1)-(B+C)z_2.   \notag
\end{align}
Let $L=Q/(B+C).$ 
Note system \eqref{4ODE} admits a globally attracting invariant line $\mathcal{L}$ on which
$$z_0=L(\al_2-\al_1), \ w_2=Ls_2(1-\al_0),  \ \mbox{ and } \ z_2=Ls_2(\al_2-\al_1)$$
are each at equilibrium. 
One can show that on $\mathcal{L}$ system \eqref{4ODE} reduces to the single equation 
\begin{equation}\label{wdot}
R\dot{w}=-B(w-F(\eta)),
\end{equation}
where we have let $w=w_0$ for ease of notation, and where
\begin{equation}\label{Fofeta}
F(\eta)=\frac{1}{B}\left(Q(1-\al_0)-A+CL(\al_2-\al_1)(\eta-\textstyle{\frac{1}{2}}+s_2P_2(\eta))\right).
\end{equation}
We note $F(\eta)$ is a cubic polynomial. Thus, equation \eqref{wdot} provides the approximation to equation \eqref{budwidA}.

For equation \eqref{budwidB}, one computes the expression for $T(\eta)$ on the invariant line $\mathcal{L}$. As detailed in \cite{dickesther}, equation \eqref{budwidB}   can be written  in the form
\begin{equation}\notag
\dot{\eta}=\rho (w-G(\eta)),
\end{equation}
where
\begin{equation}\label{Gofeta}G(\eta)=-Ls_2(1-\al_0)p_2(\eta)+T_c,
\end{equation}
when restricting to $\mathcal{L}$.
Hence the infinite-dimensional system \eqref{budwid} is approximated by the system of ODEs
\begin{align}\label{weta}
\dot{w}&=-\tau(w-F(\eta))\\\notag
\dot{\eta}&=\rho(w-G(\eta)),\notag
\end{align}
with $F(\eta)$ given by \eqref{Fofeta}, $G(\eta)$ the quadratic polynomial \eqref{Gofeta}, and $\tau=B/R$. One can show that for fixed $\eta$ the variable $w$ is a translate of the global average temperature \cite{dickesther}.

Phase planes for system \eqref{weta} are shown in Figure 1. Note the existence of equilibrium solutions is independent of $\rho$ and   $\tau$. Moreover, in \cite{dickesther} it is shown that there exist a stable equilibrium point with ice line near the pole, and a saddle equilibrium point with a large ice cap, for all $\rho>0$ (as in Figure 1(a)), for standard parameter values (also see $\S5.1$). This simple two-dimensional model captures the essence of the dynamics of the infinite-dimensional Budyko-Widiasih model \eqref{budwid}.

\begin{figure}[t!]
\includegraphics[width=5.6in,trim = 1.0in 2.5in .7in  1.6in, clip]{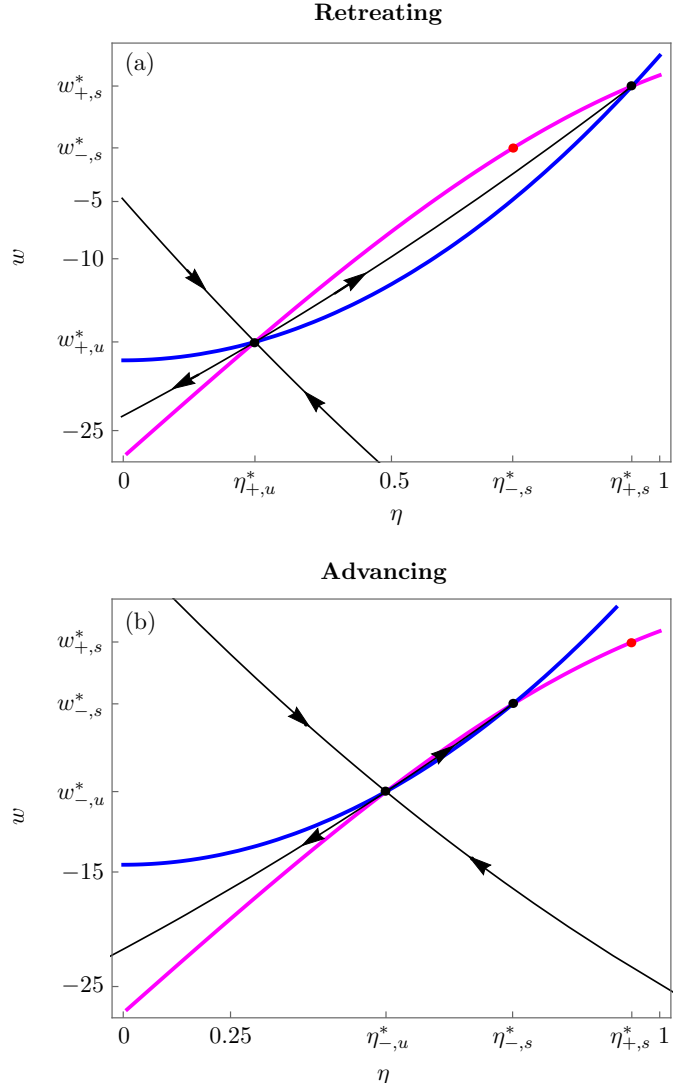}
\begin{center}
\caption{  {\small (a) Phase plane for the subsystem \eqref{retA}-\eqref{retB}, for which $T^+_c=-10^\circ$C. \ (b) Phase plane for the subsystem \eqref{advA}-\eqref{advB}, for which $T^-_c=-5.5^\circ$C. \ 
In each plot the magenta curve is the $w$-nullcline, and the blue curve is the $\eta$-nullcline.  The $\eta$-nullcline in (b) is the $\eta$-nullcline in (a) translated vertically by $T^-_c-T^+_c$. The $w$-nullclines in (a) and (b) are identical.
} }
 \end{center}
\end{figure}

As mentioned above, system \eqref{weta} does not allow for   glacial cycles, given that $\eta$ either  approaches an equilibrium position or the equator over time. Motivated by the role played by the relative sizes of the accumulation and ablation zones in glacial advance and retreat, we introduce a new snow line variable, independent  of the edge of the ice sheet, thereby supplementing equations \eqref{weta}. We present this model in the following section.

\section{Addition of a snow line}

The accumulation and ablation of ice play a  fundamental role in the theory of glacial cycles, serving to control  the terminus advance and retreat, the ice volume, and the geometry of the surface of the ice sheet \cite{bahr}. Significant motivation for our model   stems from the study \cite{abe}, in which numerical simulations of the glacial cycles over the past 400 kyr were carried out using an ice sheet model coupled to a general circulation model. These simulations typically produced two stable equilibrium ice line positions for a given set of parameters, one small and the other much larger. In addition, the larger the ice sheet at equilibrium, the larger was the ablation zone,
leading to increased ice sheet instability with regard to Milankovitch forcing. 

Most importantly, however, Abe-Ouchi et al found the fast retreat of the ice sheet was due to significantly enhanced ablation (due in turn to delayed isostatic rebound \cite{abe}). That is, the ablation rate for a large, advancing ice sheet was necessarily much smaller than the ablation rate for a retreating ice sheet, in order to faithfully reproduce the last four glacial cycles. We incorporate this simple idea into our conceptual model below. For further insight into the important role played by ablation rates (and the ratio of ablation rates to accumulation rates) in modeling glacial cycles see, for example, \cite{ghil}, \cite{kallen}, \cite{oer}, \cite{rudd}, \cite{tzip2003}, \cite{weert}.

Motivation for our model was also drawn from a  model of the thermohaline circulation presented in \cite{welander}. This simple ``flip-flop" conceptual model exhibits oscillations in the temperature and salinity of a well-mixed ocean layer. The model vector field has a line of discontinuity
 that produces a switch to the alternate regime when intersected by a trajectory. Similarities between our model and that presented in \cite{welander} will become evident below.

Consider once again system \eqref{weta}. Recall that for fixed $\eta$-values, $w$ is a translate of the global average temperature. We introduce independent snow and ice lines so as to incorporate accumulation and ablation zones. The model continues to depend upon latitude as discussed in Sections 2 and 3.

We begin by recasting the role played by $\eta$, interpreting $\eta$ henceforth as the snow line. We let $\xi$ denote the ice line, that is, the edge of the ice sheet (see Figure 2). The ablation zone has extent $\eta-\xi$ (when $\eta>\xi$), while the accumulation zone has size $1-\eta$. 

\vspace*{0.5in}
\begin{figure}[h!]
\begin{center}
\includegraphics[width=5.25in,trim = 1in 7.6in 1in  .9in, clip]{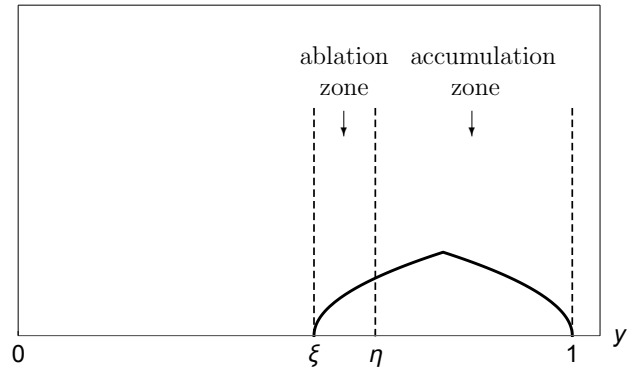}\\
\caption{{\small The model set-up. \ $\eta$ is the snow line and $\xi$ is the ice line. The shape of the glacier is for illustrative purposes only.}} 
 \end{center}
\end{figure}

We continue to use albedo function \eqref{budalb}. We note the albedo of the portion of the ice sheet in the ablation zone is much smaller than that of the ice sheet lying poleward of $\eta$. For example,  as the old ice emerges and melts in the ablation zone, a surface layer of dust that was originally laid down with snowfall high on the ice sheet
is created \cite{bog}. Additionally, locally increased winds further cover the surface in the ablation zone with dust and debris \cite{pelt}. Albedo value $\al_1$ in \eqref{budalb} might then be interpreted as an average  of the albedos of the planet's surface south of $\xi$ and the darkened (relative to snow) ice surface in the ablation zone.

The use of   independent ice and snow lines was also partially due to the desire to include ice volume in this latitude-dependent model. (A strong argument for the inclusion of ice volume in any EBM is given in \cite{oer}.) While our model does not explicitly include volume, the latter can be recovered from $\xi$  were we to assume a specific ice sheet shape as, for example,  in \cite{weert}. Finally, the model allows for different time scales for the snow line and the (more slowly moving) ice line. 

\subsection{Model equations}

The temperature-ice line-snow line model is a nonsmooth system with state space
$$\mathcal{B}=\{(w,\eta,\xi): w\in\R, \eta\in [0,1], \xi\in [0,1]\},$$defined as  follows. Pick parameters $b_0<b<b_1$   representing  ablation rates, and a parameter $a$ denoting the accumulation rate. When $b(\eta-\xi)-a(1-\eta)<0$, so that accumulation exceeds ablation and the ice sheet advances, set
\begin{subequations}\label{adv}
\begin{align}
\dot{w}&=-\tau(w-F(\eta))=f^-_1(w,\eta,\xi)  \label{advA} \\
\dot{\eta}&=\rho(w-G_-(\eta))=f^-_2(w,\eta,\xi)   \label{advB}  \\
\dot{\xi}&=\eps(b_0(\eta-\xi)-a(1-\eta))=f^-_3(w,\eta,\xi).  \label{advC}
\end{align}
\end{subequations}
The function $F(\eta)$ in \eqref{advA} is given by  \eqref{Fofeta},   $G_-(\eta)$ in \eqref{advB} is given by  \eqref{Gofeta}  but  with $T_c=T^-_c=-5.5^\circ$C, and $\eps>0$ is a time constant for the movement of the ice line.

When ablation exceeds accumulation ($b(\eta-\xi)-a(1-\eta)>0$) and the ice sheet retreats, we set
\begin{subequations}\label{ret}
\begin{align}
\dot{w}&=-\tau(w-F(\eta))=f^+_1(w,\eta,\xi)  \label{retA} \\
\dot{\eta}&=\rho(w-G_+(\eta))=f^+_2(w,\eta,\xi)   \label{retB}  \\
\dot{\xi}&=\eps(b_1(\eta-\xi)-a(1-\eta))=f^+_3(w,\eta,\xi).  \label{retC}
\end{align}
\end{subequations}
The function $F(\eta)$ in \eqref{retA} is again given by  \eqref{Fofeta}, while $G_+(\eta)$ in \eqref{retB} is given by  \eqref{Gofeta}   with $T_c=T^+_c=-10^\circ$C.  

The relative sizes of ablation rates $b_0$ and $b_1$ were motivated by \cite{abe} as described above. The choice of different  $T_c$-values  is motivated by \cite{tzip2003}, in which a linear interpolation  between $T_c=-13^\circ$C and $T_c=-3^\circ$C is introduced to model changes in deep ocean temperature. The intuitive idea is that a large advancing ice sheet implies a colder world generally, so that less energy is require to form ice (and vice versa for a retreating ice sheet).

We thus arrive at a three-dimensional system having   a plane of discontinuity
\begin{equation}\label{Sig}
\Sigma=\{ \, (w,\eta,\xi) : b(\eta-\xi)-a(1-\eta)=0 \, \} = \{ \, (w,\eta,\xi) : \xi=(1+\textstyle{\frac{a}{b}})\eta-\textstyle{\frac{a}{b}}\equiv \gamma(\eta) \, \}.
\end{equation}
As we will see, a trajectory in $(w,\eta,\xi)$-space passing through $\Sigma$ switches from advancing mode to glacial retreat, or vice versa, similar in spirit to the flip-flop model in \cite{welander}.

Let  $$S_+=\{ \, (w,\eta,\xi)\in \mathcal{B} :  \xi<\gamma(\eta) \, \} \mbox{ \ and \ } S_-=\{ \, (w,\eta,\xi)\in\mathcal{B} :  \xi>\gamma(\eta) \, \},$$
with $\gamma(\eta)$   as in \eqref{Sig}. Called a  {\em discontinuity boundary} (or {\em switching boundary}), $\Sigma$ is a hyperplane  separating $\mathcal{B}$ into subspaces $S_+$ and $S_-$. Consider the vector fields
$${\bf V}_\pm:S_\pm\rightarrow \R^3, \ {\bf V}_\pm(w,\eta,\xi)=(f^\pm_1(w,\eta,\xi), f^\pm_2(w,\eta,\xi),f^\pm_3(w,\eta,\xi)),$$
noting ${\bf V}_\pm$ each extend smoothly to $\Sigma$. For ${\bf x}=(w,\eta,\xi)\in \mathcal{B},$ consider the differential inclusion
\begin{equation}\label{Fil}
\dot{\bf x}\in {\bf V}({\bf x})=
\begin{cases}
{\bf V}_-({\bf x}), &  {\bf x}\in S_- \\ 
\{ (1-q){\bf V}_-({\bf x})+q {\bf V}_+({\bf x}) : q\in [0,1]\}, & {\bf x}\in\Sigma \\ 
{\bf V}_+({\bf x}), &  {\bf x}\in S_+ . 
\end{cases}
\end{equation}
Solutions, while in $S_-$, have flow $\phi_-({\bf x},t)$ given by $\dot{\bf x}={\bf V}_-({\bf x})$, and solutions have flow $\phi_+({\bf x},t)$
given by  $\dot{\bf x}={\bf V}_+({\bf x})$ while in $S_+$.  For ${\bf x}\in\Sigma, \ \dot{\bf x}$ must lie in the closed convex hull of the two vectors ${\bf V}_-({\bf x})$ and ${\bf V}_+({\bf x})$.

A solution to \eqref{Fil} {\em in the sense of Filippov} is an absolutely continuous function ${\bf x}(t)$ satisfying $\dot{\bf x}\in {\bf V}({\bf x})$ for almost all $t$. (Note $\dot{\bf x}(t)$ is not defined at times for which ${\bf x}(t)$ arrives at or leaves $\Sigma$.) Given that ${\bf V}_\pm$ are continuous on $S_\pm\cup \Sigma$, the set-valued map ${\bf V}({\bf x})$ is upper semi-continuous, and closed, convex and bounded for all ${\bf x}\in \mathcal{B}$ and $t\in\R$. This implies that for each ${\bf x}_0\in \mbox{Int}(\mathcal{B})$ there is a solution ${\bf x}(t)$ to differential inclusion \eqref{Fil} in the sense of Filippov, defined  on an interval $[0,t_f]$, with ${\bf x}(0)={\bf x}_0$   \cite{leine}.

We now turn to an analysis of the model.

\section{Model analysis}

We first note   model equations \eqref{adv} and \eqref{ret} partially decouple: As long as a trajectory remains in $S_-$ or $S_+$, changes in $\xi$ play no role in the evolution of $w$ and $\eta$. This simple observation will prove useful in the following analysis.

We set several parameter values before continuing.

\vspace{0.2in}
\begin{center}{\small
\begin{tabular}{|ccc|ccc|}\hline
Parameter & Value & Units & Parameter & Value & Units\\ \hline
$Q$ & 343 & Wm$^{-2}$ & $T^-_c$ & $-5.5$ & $^\circ$C\\ 
$A$ & 202 & Wm$^{-2}$ & $b_0$ & 1.5 & dimensionless\\ 
$B$ & 1.9 & Wm$^{-2}(^\circ\mbox{C})^{-1}$ & $b$ & 1.75 & dimensionless\\
$C$ & 3.04 & Wm$^{-2}(^\circ\mbox{C})^{-1}$ & $b_1$ & 5 & dimensionless\\
$\al_1$ & 0.32 & dimensionless & $a$ & 1.05 & dimensionless\\
$\al_2$ & 0.62 & dimensionless & $\tau$ & 1 & seconds$^{-1}$\\ 
$T^+_c$ & $-10$ & $^\circ$C & $\rho$ & 0.1 & seconds$^{-1} (^\circ\mbox{C})^{-1}$\\ \hline
\end{tabular}

\vspace{0.2in}
{\bf Table 1} Parameter values (left column from \cite{dickesther}). The parameter $\eps$ has units (seconds)$^{-1}$.
}
\end{center}

\vspace{0.1in}
\subsection{Equilibria}

Note time constants $\tau, \rho$ and $\eps$ play no role in the determination of equilibrium solutions of systems \eqref{adv} and \eqref{ret}. As we'll see below, they also play no role in determining the stability type of any equilibria. 

An equilibrium point ${\bf Q}^*_+=(w^*_+,\eta^*_+,\xi^*_+)$ for \eqref{ret} satisfies $w^*_+=F(\eta^*_+)=G_+(\eta^*_+)$ (that is, $F(\eta^*_+)-G_+(\eta^*_+)=0$), and $\xi^*_+=(1+\frac{a}{b_1})\eta^*_+-\frac{a}{b_1}$. In \cite{dickesther} it was shown that, with parameters as in Table 1, the (decoupled) subsystem \eqref{retA}-\eqref{retB} has a ``small ice cap" sink $(w^*_{+,s},\eta^*_{+,s})\approx (5.08, 0.95)$ and a ``large ice cap" saddle $(w^*_{+,u},\eta^*_{+,u})\approx (-17.26, 0.25)$ (see Figure 1(a)), and no other equilibria having  $\eta\in[0,1]$. Hence system \eqref{ret} has two equilibria
$${\bf Q}^*_{+,s}=(w^*_{+,s},\eta^*_{+,s}, (1+\textstyle{\frac{a}{b_1}})\eta^*_{+,s}-\textstyle{\frac{a}{b_1}}) \ \mbox{ and } \ 
{\bf Q}^*_{+,u}=(w^*_{+,u},\eta^*_{+,u}, (1+\textstyle{\frac{a}{b_1}})\eta^*_{+,u}-\textstyle{\frac{a}{b_1}}).
$$
The Jacobian matrix for \eqref{ret} is
\begin{equation}\label{Jret}
J_+(w,\eta,\xi)=\left[\begin{matrix} -\tau& \tau F^\prime(\eta) & 0\\ \rho& -\rho G^\prime_+(\eta) & 0\\ 0 & \eps(b_1+a)& -\eps b_1\end{matrix}\right].
\end{equation}
Thus at any equilibrium point ${\bf Q}^*$, \ $-\eps b_1$ is an eigenvalue of $J_+({\bf Q}^*)$ with eigenvector $[0 \ 0\ 1]^T$. Combining this with the analysis in \cite{dickesther}, we see that $J_+({\bf Q}^*_{+,s})$ has three negative eigenvalues, while $J_+({\bf Q}^*_{+,u})$
has two negative eigenvalues and one positive eigenvalue. We have that ${\bf Q}^*_{+,s}$ is a sink, and ${\bf Q}^*_{+,u}$ is a saddle with two-dimensional stable manifold, for system \eqref{ret}.

We note
$$\xi^*_{+,s}=(1+\textstyle{\frac{a}{b_1}})\eta^*_{+,s}-\textstyle{\frac{a}{b_1}}>(1+\textstyle{\frac{a}{b}})\eta^*_{+,s}-\textstyle{\frac{a}{b}}=\gamma(\eta^*_{+,s}),$$
implying that ${\bf Q}^*_{+,s}\in S_-$. Such an equilibrium point is called a  virtual equilibrium point \cite{dibernardo}, in the sense that ${\bf Q}^*_{+,s}$ cannot be ``seen" by the vector field ${\bf V}_+$. A $\phi_+$-trajectory in $S_+$ will intersect $\Sigma$ (and switch to the flow $\phi_-$) before having the chance to approach ${\bf Q}^*_{+,s}$.

\vspace{0.2in}
\noindent
{\bf Definition 5.1} \ Let ${\bf x}^*\in \mathcal{B}$.\\
(i) \  ${\bf x^*}$ is a {\em regular equilibrium point} of \eqref{Fil} if either ${\bf V}_+({\bf x}^*)={\bf 0}$ and ${\bf x}^*\in S_+$, or if ${\bf V}_-({\bf x}^*)={\bf 0}$ and ${\bf x}^*\in S_-$.\\
(ii) \ ${\bf x^*}$ is a {\em virtual equilibrium point} of \eqref{Fil} if either ${\bf V}_+({\bf x}^*)={\bf 0}$ and ${\bf x}^*\in S_-$, or if ${\bf V}_-({\bf x}^*)={\bf 0}$ and ${\bf x}^*\in S_+$.\\
(iii) \ ${\bf x^*}$ is a {\em boundary equilibrium point} of \eqref{Fil} if  ${\bf V}_+({\bf x}^*)={\bf V}_-({\bf x}^*)={\bf 0}$ and ${\bf x}^*\in \Sigma$.

\vspace{0.2in}
The equilibrium point analysis for vector field ${\bf V}_-$ is similar. Note the Jacobian matrix $J_-(w,\eta,\xi)$ for \eqref{adv} is identical to \eqref{Jret}, except that $b_1$ is replaced with $b_0$ (observe $G^\prime_+(\eta)=G^\prime_-(\eta)$). As in \cite{dickesther}, the decoupled $(w,\eta)$-block
$$M=\left[\begin{matrix} -\tau& \tau F^\prime(\eta)\\ \rho & -\rho G^\prime_-(\eta)\end{matrix}\right]$$
satisfies Tr$(M)=-(\tau+\rho G^\prime_-(\eta))<0$ for $\eta>0$, while det$(M)=-\tau\rho(F^\prime(\eta)-G^\prime_-(\eta))$. One can   rigorously show the function $h_-(\eta) =F(\eta)-G_-(\eta)$ has two zeros $\eta^*_{-,u}< \eta^*_{-,s}$ in [0,1], with $h^\prime_-(\eta^*_{-,u})>0$ and $h^\prime_-(\eta^*_{-,s})<0$ (also see Figure 3). Thus subsystem \eqref{advA}-\eqref{advB} has a saddle point $(w^*_{-,u},\eta^*_{-,u})$ and a sink $(w^*_{-,s},\eta^*_{-,s})$ (see Figure 1(b)), from which it follows system \eqref{adv} admits a saddle point ${\bf Q}^*_{-,u}=(w^*_{-,u},\eta^*_{-,u}, \xi^*_{-,u})$ with two-dimensional stable manifold, and a sink ${\bf Q}^*_{-,s}=(w^*_{-,s},\eta^*_{-,s}, \xi^*_{-,s})$. Once again the sink ${\bf Q}^*_{-,s}$ is virtual as $$ \xi^*_{-,s}=(1+\textstyle{\frac{a}{b_0}})\eta^*_{-,s}-\textstyle{\frac{a}{b_0}}<(1+\textstyle{\frac{a}{b}})\eta^*_{-,s}-\textstyle{\frac{a}{b}}=\gamma(\eta^*_{-,s}),$$ implying ${\bf Q}^*_{-,s}\in S_+.$

A key observation can be gleaned from Figure 1. First note subsystems \eqref{advA}-\eqref{advB} and \eqref{retA}-\eqref{retB} differ only in that $G_-(\eta)=G_+(\eta)+(T^-_c-T^+_c)$. As indicated in Figure 3, this simply moves $\eta^*_{-,u}$ and $\eta^*_{-,s}$ closer together, relative to $\eta^*_{+,u}$ and $\eta^*_{+,s}$.
Also note   $(w^*_{-,s},\eta^*_{-,s})$ is in the stable set of $(w^*_{+,s},\eta^*_{+,s})$ under the two-dimensional  flow corresponding to subsystem \eqref{retA}-\eqref{retB}, while  $(w^*_{+,s},\eta^*_{+,s})$ is in the stable set of $(w^*_{-,s},\eta^*_{-,s})$ under the two-dimensional  flow corresponding to subsystem \eqref{advA}-\eqref{advB}. Given the linear nature of equations \eqref{advC}-\eqref{retC}  (along with the aforementioned decoupling), it follows that   ${\bf Q}^*_{-,s}$ is in the stable set of ${\bf Q}^*_{+,s}$ under the flow $\phi_+$  (system \eqref{ret}), and ${\bf Q}^*_{+,s}$ is in the stable set of ${\bf Q}^*_{-,s}$ under the flow $\phi_-$  (system \eqref{adv}).

More generally, the stable set of $(w^*_{-,s},\eta^*_{-,s})$ under subsystem \eqref{advA}-\eqref{advB} is a subset of  the stable set of $(w^*_{+,s},\eta^*_{+,s})$ under subsystem \eqref{retA}-\eqref{retB}. Thus the stable set of ${\bf Q}^*_{-,s}$ under the flow $\phi_-$ (system \eqref{adv}) is a subset of the stable set of ${\bf Q}^*_{+,s}$ under the flow $\phi_+$ (system \eqref{ret}).

We now turn to the behavior of trajectories intersecting the discontinuity boundary $\Sigma$.

\begin{figure}[t!]
\begin{center}
\includegraphics[width=5.7in,trim = 1in 6.43in 1in  1.7in, clip]{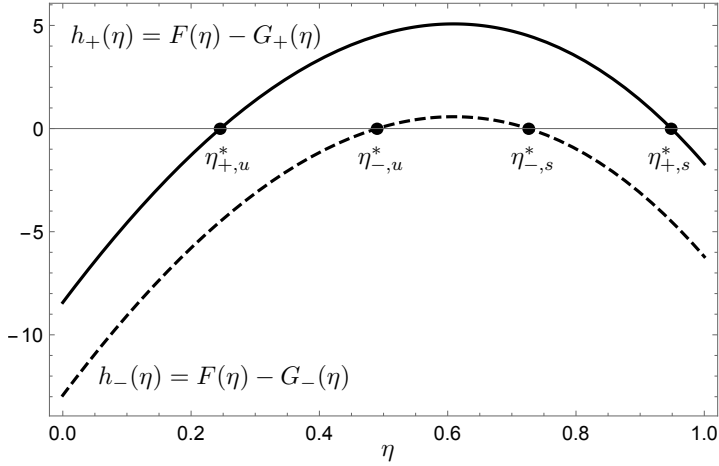}\\
\caption{{\small The $\eta$-coordinates at equilibrium. {\em Solid: } Subsystem \eqref{retA}-\eqref{retB}, for which $T^+_c=-10^\circ$C. \ {\em Dashed: }  Subsystem \eqref{advA}-\eqref{advB}, for which $T^-_c=-5.5^\circ$C.}}
 \end{center}
 
 \vspace*{0.2in}
\end{figure}

\subsection{Behavior on the discontinuity boundary $\boldsymbol{\Sigma}$}

The   discontinuity boundary $\Sigma$ given by  \eqref{Sig} is a plane with normal vector $${\bf N}=\left[\begin{matrix} w & \eta & \xi\end{matrix}\right]^T=
\left[\begin{matrix} 0 &1+\textstyle{\frac{a}{b}} & -1\end{matrix}\right]^T.$$
For ${\bf x}\in \Sigma$, a computation yields ${\bf V}_+({\bf x}) \perp {\bf N}$ if and only if 
\begin{equation}\label{gplus}
w=G_+(\eta)+\frac{\eps a(1-\eta)(b_1-b)}{\rho(a+b)}\equiv g_+(\eta).
\end{equation}
Thus ${\bf V}_+$ is tangent to $\Sigma$ only on the parabola \ 
$\Lambda_+=\{(w,\eta,\xi)\in\Sigma : w=g_+(\eta)\}$ (the red curve in Figure 4). Furthermore, given ${\bf x}=(w,\eta,\xi)\in\Sigma, \ {\bf V}_+({\bf x}) \dotp {\bf N}<0$ for $w<g_+(\eta)$, and  ${\bf V}_+({\bf x}) \dotp {\bf N}>0$ for $w>g_+(\eta)$. Referring to the orientation provided in Figure 4, $\phi_+$ trajectories pass through $\Sigma$ below $\Lambda_+$ from front to back, while $\phi_+$ trajectories pass through $\Sigma$ above $\Lambda_+$ from back to front.

Similar computations reveal that ${\bf V}_-$ is tangent to $\Sigma$ if and only if
\begin{equation}\label{gminus}
w=G_-(\eta)+\frac{\eps a(1-\eta)(b_0-b)}{\rho(a+b)}\equiv g_-(\eta),
\end{equation}
that is, only on the parabola  \ $\Lambda_-=\{(w,\eta,\xi)\in\Sigma : w=g_-(\eta)\}$ (the blue curve in Figure 4). 
 Given ${\bf x}=(w,\eta,\xi)\in\Sigma, \ {\bf V}_-({\bf x}) \dotp {\bf N}>0$ for $w>g_-(\eta)$, and  ${\bf V}_+({\bf x}) \dotp {\bf N}<0$ for $w<g_+(\eta)$. Thus $\phi_-$ trajectories pass through $\Sigma$ below $\Lambda_-$ from front to back, while $\phi_-$ trajectories pass through $\Sigma$ above $\Lambda_-$ from back to front.

Let $$\Sigma^{\mbox{\scriptsize SL}}=\{  (w,\eta,\xi) \in \Sigma: g_+(\eta)<w<g_-(\eta)  \}.$$ For ${\bf x}\in \Sigma^{\mbox{\scriptsize SL}}, {\bf V}_+({\bf x})\dotp {\bf N}>0$ (so ${\bf V}_+({\bf x})$ points into $S_+$), while ${\bf V}_-({\bf x})\dotp {\bf N}<0$ (so ${\bf V}_-({\bf x})$ points into $S_-$). The subset $\Sigma^{\mbox{\scriptsize SL}}$ of the discontinuity boundary $\Sigma$ is therefore a {\em  repelling sliding region} \cite{leine}; Filippov's approach does not provide for unique solutions ${\bf x}(t)$ in forward time if ${\bf x}(0)\in 
\Sigma^{\mbox{\scriptsize SL}}$ \cite{fil}.

Let $$\Sigma_+=\{  (w,\eta,\xi) \in \Sigma: w<g_+(\eta)  \} \mbox{ \ and \ } \Sigma_-=\{  (w,\eta,\xi) \in \Sigma: w>g_-(\eta)  \}.$$ Unique solutions to \eqref{Fil} do exist for trajectories passing through $\Sigma_+$ or $\Sigma_-$, as ${\bf V}_-({\bf x})\dotp {\bf N}$ and ${\bf V}_+({\bf x})\dotp {\bf N}$ have the same sign on each of $\Sigma_+$ and $\Sigma_-$.
 A $\phi_+$-trajectory in $S_+$ intersecting $\Sigma_+$ at ${\bf x}$ will cross $\Sigma$ transversally, becoming a $\phi_-$-trajectory at ${\bf x}$. Similarly, a $\phi_-$-trajectory in $S_-$ intersecting $\Sigma_-$ at ${\bf y}$ will cross $\Sigma$ transversally, becoming a $\phi_+$-trajectory at ${\bf y}$. In this scenario a unique (though nonsmooth) solution of \eqref{Fil} exists, in the sense of Filippov.

\subsection{Section maps for the Filippov flow}

\begin{figure}[t!]
\begin{center}
\includegraphics[width=6in,trim = 1in 4.8in 1in  1.9in, clip]{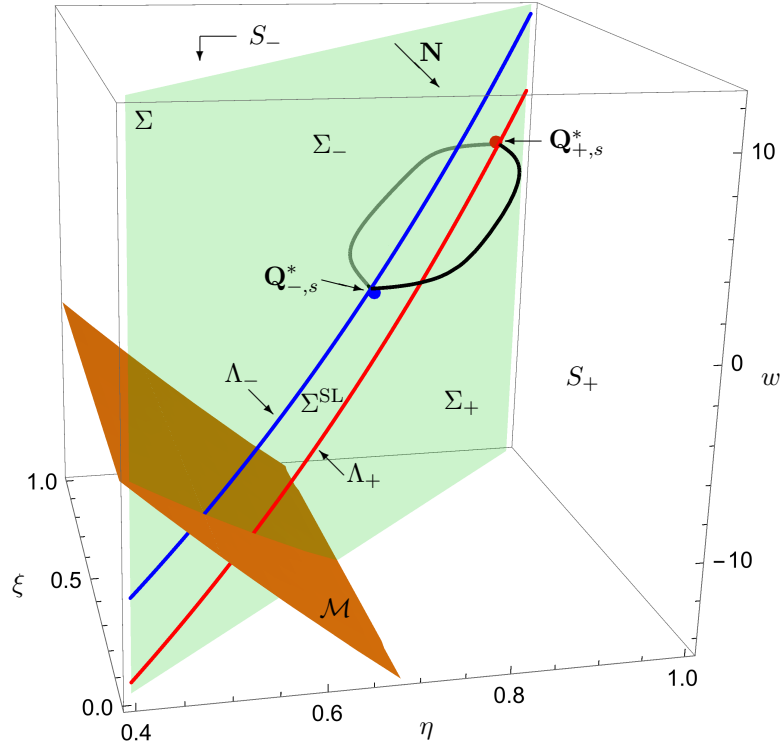}\\
\caption{{\small The discontinuity boundary  $\Sigma$ and domains $S_\pm$ for vector fields ${\bf V}_\pm$. ${\bf Q}^*_{-,s}\in S_+$ and ${\bf Q}^*_{+,s}\in S_-$  are virtual equilibria. Also pictured is the periodic orbit whose existence follows from Theorem 5.6.}}
 \end{center}
 
 \vspace*{0.2in}
\end{figure}

One can show the curves of tangency $\Lambda_\pm$ of ${\bf V}_\pm$ with $\Sigma$ intersect precisely when 
\begin{equation}\label{etaeps}
\eta=\eta(\eps)=1-\frac{(T^-_c-T^+_c)\rho(a+b)}{\eps a(b_1-b_0)}.
\end{equation}
Thus if
\begin{equation}\label{cond1}
\eps<\frac{(T^-_c-T^+_c)\rho(a+b)}{ a(b_1-b_0)},
\end{equation}
$\eta(\eps)$ in \eqref{etaeps} is strictly negative. We assume $\eps$ satisfies this bound in all that follows, ensuring $\Lambda_+$ and $\Lambda_-$ do not intersect in $\mathcal{B}$ (as in Figure 4).

Plotted in Figure 4 is the two-dimensional stable manifold ${\cal M}$ of the saddle point ${\bf Q}^*_{-,u}$ under the flow $\phi_-$. Due to the $(w,\eta)$-decoupling mentioned previously,
the projection of ${\cal M}$ onto the $(w,\eta)$-plane is simply the (one-dimensional) stable manifold of $(w^*_{-,u},\eta^*_{-,u})$ under the flow given by   subsystem \eqref{advA}-\eqref{advB} and plotted in Figure 1(b). In particular, we note the plot of the surface ${\cal M}$ is independent of $\eps>0$.

${\cal M}$ partitions $\Sigma$ into two subsets, a ``top" (larger $w$) and a ``bottom" (smaller $w$). Denote this top subset by $\mathcal{T}$. For any point ${\bf x}=(w,\eta,\xi)\in \mathcal{T}, \phi_-({\bf x},t)=(w(t),\eta(t),\xi(t))\rightarrow {\bf Q}^*_{-,s}$. This follows from the fact $(w(t),\eta(t))\rightarrow (w^*_{-,s},\eta^*_{-,s})$ as $t\rightarrow\infty$, and hence $\xi(t)\rightarrow \xi^*_{-,s}$ by equation \eqref{advC}. Thus $\mathcal{T}\subset W^s_-({\bf Q}^*_{-,s})$, the stable set of the sink ${\bf Q}^*_{-,s}$ under the flow $\phi_-$. By the observations made in $\S$5.1, it follows $\mathcal{T}\subset W^s_+({\bf Q}^*_{+,s})$, the stable set of the sink ${\bf Q}^*_{+,s}$ under the flow $\phi_+$. Again, all of this holds for any $\eps>0$.

We have that for any $\eps>0, \ \Sigma_+\cap \mathcal{T}\subset W^s_-({\bf Q}^*_{-,s})$  and ${\bf V}_-({\bf x}) \dotp {\bf N}<0$ for all ${\bf x}\in\Sigma_+\cap \mathcal{T}$. Hence for any $\eps>0$,  the $\phi_-$-trajectory starting at ${\bf x}\in \Sigma_+\cap \mathcal{T}$ enters $S_-$ before necessarily intersecting  $\Sigma_-\cap \mathcal{T}$ as it seeks to approach the (virtual) sink ${\bf Q}^*_{-,s}$.  For any $\eps>0$ and for each ${\bf x}\in \Sigma_+\cap \mathcal{T}$,  there is then a $t=t({\bf x})>0$ such that $\phi_-({\bf x},t)\in\Sigma_-\cap\mathcal{T}$.  We write $t({\bf x})=t({\bf x},\eps)$ to emphasize the dependence of this $t$-value on $\eps$ as well as on ${\bf x}$. We note, for future reference, that $t({\bf x},\eps)\rightarrow\infty $ as $\eps\rightarrow 0^+.$

This provides, for any $\eps>0$, a continuous   map $$r^-_\eps:\Sigma_+\cap \mathcal{T}\rightarrow \Sigma_-\cap \mathcal{T}, \ r^-_\eps({\bf x})=\phi_-({\bf x}, t({\bf x}, \eps)).$$

Similarly,  for any $\eps>0, \ \Sigma_-\cap \mathcal{T}\subset W^s_+({\bf Q}^*_{+,s})$ and ${\bf V}_+({\bf y}) \dotp {\bf N}>0$ for all ${\bf y}\in \Sigma_-\cap \mathcal{T}.$ Hence for any $\eps>0$,  the $\phi_+$-trajectory starting at ${\bf y}\in \Sigma_-\cap \mathcal{T}$ enters $S_+$ before necessarily intersecting  $\Sigma_+\cap \mathcal{T}$ as it tries to approach the (virtual) sink ${\bf Q}^*_{+,s}$.  For any $\eps>0$ and for each ${\bf y}\in \Sigma_-\cap \mathcal{T}$,  there exists a $t=t({\bf y},\eps)>0$ such that $\phi_+({\bf y},t)\in\Sigma_+\cap\mathcal{T}$ (again,   $t({\bf y},\eps)\rightarrow\infty $ as $\eps\rightarrow 0^+$). 
This defines 
a continuous   map $$r^+_\eps:\Sigma_-\cap \mathcal{T}\rightarrow \Sigma_+\cap\mathcal{T}, \ r^+_\eps({\bf y})=\phi_+({\bf y}, t({\bf y}, \eps)),$$
which  exists for any $\eps>0$. 

We then have two continuous section maps
\begin{equation}\label{bothways}
r^+_\eps\circ r^-_\eps : \Sigma_+\cap \mathcal{T} \rightarrow \Sigma_+\cap \mathcal{T} \mbox{  \ and  \ }
r^-_\eps\circ r^+_\eps : \Sigma_-\cap \mathcal{T} \rightarrow  \Sigma_-\cap \mathcal{T},
\end{equation}
each defined for any $\eps>0$. Given certain hypotheses, these maps are each contraction maps, as we show in the following section.

\subsection{Existence of an attracting periodic orbit}

We begin with the following propositions.

\vspace{0.2in}
\noindent
{\bf Proposition 5.4} \ (i) Given $c_1\in (0,1)$ and a compact set $D_+\subset \Sigma_+\cap \mathcal{T}$, there exists $\eps_1>0$ such that for all $\eps\leq \eps_1$ and for all ${\bf x}_1, {\bf x}_2\in D_+,$
 \ $\|r^-_\eps({\bf x}_2)-r^-_\eps({\bf x}_1)\|\leq c_1\|{\bf x}_2-{\bf x}_1\|.$
 
 \vspace{0.1in}
 \noindent
 (ii) Given $c_2\in (0,1)$ and a compact set $D_-\subset \Sigma_-\cap \mathcal{T}$, there exists $\eps_2>0$ such that for all $\eps\leq \eps_2$ and for all ${\bf y}_1, {\bf y}_2\in D_-,$
 \ $\|r^+_\eps({\bf y}_2)-r^+_\eps({\bf y}_1)\|\leq c_2\|{\bf y}_2-{\bf y}_1\|.$

\vspace{0.2in}
\noindent
{\em Proof.} We prove (i). Let $c_1\in (0,1)$ and let $D_+\subset \Sigma_+\cap \mathcal{T}, \ D_+$ compact. Let $\psi_-$ denote the flow corresponding to subsystem \eqref{advA}-\eqref{advB} (with phase plane depicted in Figure 1(b)). For any ${\bf x}=(w_0,\eta_0,\xi_0)\in D_+\subset\mathcal{T}, \ \psi_-((w_0,\eta_0),t)\rightarrow (w^*_{-,s},\eta^*_{-,s})$ as $t\rightarrow\infty$. As $E=\{(w,\eta) : (w,\eta,\xi)\in D_+\}$ is compact, there exists $T_1$ such that for all $t\geq T_1$ and for all ${\bf u}, {\bf v}\in E$,
\begin{equation}\label{psicontract}
\|\psi_-({\bf u},t)-\psi_-({\bf v},t)\|\leq c_1\|{\bf u}-{\bf v}\|.
\end{equation}
Given ${\bf x}\in D_+$, pick $\eps({\bf x})>0$ such that $t({\bf x},\eps({\bf x}))>T_1$ (recalling $\phi_-({\bf x},t({\bf x},\eps({\bf x})))  
=r^-_{\eps({\bf x})}({\bf x})\in \Sigma_-\cap\mathcal{T}).$ 

By the continuity of $\phi_-$ with respect to initial conditions and time, there exists $\delta({\bf x})>0$ so that for all ${\bf y}\in B_{\delta({\bf x})}({\bf x}), \ t({\bf y},\eps({\bf x}))>T_1$, where $r^-_{\eps({\bf x})}({\bf y})\in \Sigma_-\cap\mathcal{T}.$ Note for all $\eps\leq \eps({\bf x}), t({\bf y},\eps)>T_1$.

Letting ${\bf x}$ vary, we get an open covering 
$$D_+\subset \bigcup_{{\bf x}\in D_+} B_{\delta({\bf x})}({\bf x})$$
of the compact set $D_+$. Let $\{B_{\delta({\bf x}_i)}({\bf x}_i) : i=1, ... , N\}$ be a finite subcover, and set $\eps_1=\min\{\eps({\bf x}_i): i=1, ... , N\}$. Then for any $\eps \leq \eps_1$ and for all ${\bf x}\in D_+, \ t({\bf x},\eps)>T_1$, where $r^-_{\eps}({\bf x})\in \Sigma_-\cap\mathcal{T}.$

Now let $\eps\leq \eps_1$, and let ${\bf x}_1=(w_1,\eta_1,\gamma(\eta_1)), {\bf x}_2=(w_2,\eta_2,\gamma(\eta_2))\in D_+$. Set ${\bf u}=(w_1,\eta_1)$ and 
${\bf v}=(w_2,\eta_2)$. Let ${\bf y}_1=r^-_\eps({\bf x_1})=(w^\prime_1,\eta^\prime_1,\gamma(\eta^\prime_1))$ and 
${\bf y}_2=r^-_\eps({\bf x_2})=(w^\prime_2,\eta^\prime_2,\gamma(\eta^\prime_2))$. Set ${\bf u}^\pr=(w^\pr_1,\eta^\pr_1)$ and 
${\bf v}^\pr=(w^\pr_2,\eta^\pr_2)$. By \eqref{psicontract} and our choice of $\eps$, $ \|{\bf v}^\pr-{\bf u}^\pr\|^2\leq c^2_1\|{\bf v}-{\bf u}\|^2$. We have 
\begin{align}\notag
\|{\bf y}_2-{\bf y}_1\|^2 & = \|{\bf v}^\pr-{\bf u}^\pr\|^2+(\gamma(\eta^\pr_2)-\gamma(\eta^\pr_1))^2\\\notag
& = \|{\bf v}^\pr-{\bf u}^\pr\|^2+ (1+\textstyle{\frac{a}{b}})^2 (\eta^\pr_2- \eta^\pr_1)^2\\\notag
 & \leq c^2_1\|{\bf v}-{\bf u}\|^2+(1+\textstyle{\frac{a}{b}})^2 c^2_1 (\eta_2-\eta_1)^2\\\notag
&= c^2_1\|{\bf v}-{\bf u}\|^2+c^2_1(\gamma(\eta_2)-\gamma(\eta_1))^2\\\notag
&=c^2_1\|{\bf x}_2-{\bf x}_1\|^2.
\end{align}
Hence statement (i) holds. The proof of statement (ii) is similar. \ $_\square$

\vspace{0.2in}
To use Proposition 5.4 to  show that, say, $r^+_\eps\circ r^-_\eps$ is a contraction map, 
we need a compact set $D_+\subset \Sigma_+\cap\mathcal{T}$ such that $r^+_\eps\circ r^-_\eps(D_+)\subset D_+$. To that end, let ${\bf Z}^*_+=(w^*_{+,s},\eta^*_{+,s},\gamma(\eta^*_{+,s}))$, and note ${\bf Z}^*_+\in \mathcal{T}\cap\Sigma=\mathcal{T}$. Moreover, that ${\bf Z}^*_+\in \Sigma_+$ follows from the fact $w^*_{+,s}=G_+(\eta^*_{+,s})<g_+(\eta^*_{+,s})$, which in turn is a consequence of \eqref{gplus}.

Similarly, let   ${\bf Z}^*_-=(w^*_{-,s},\eta^*_{-,s},\gamma(\eta^*_{-,s}))\in \mathcal{T} $.   As $w^*_{-,s}=G_-(\eta^*_{-,s})>g_-(\eta^*_{-,s})$ by equation \eqref{gminus},  we have ${\bf Z}^*_-\in \Sigma_-$.

\vspace{0.2in}
\noindent
{\bf Proposition 5.5} \ (i) Let $D_+$ be any compact subset of $\Sigma_+\cap\mathcal{T}$ with ${\bf Z}^*_+\in \mbox{Int}(D_+)$. Let $D_-$ be any compact subset of $\Sigma_-\cap\mathcal{T}$. There exists $\eps_3>0$ such that for all $\eps\leq \eps_3, \ r^+_\eps(D_-)\subset D_+.$

\vspace{0.051in}
\noindent
(ii) Let $D_-$ be any compact subset of $\Sigma_-\cap\mathcal{T}$ with ${\bf Z}^*_-\in \mbox{Int}(D_-)$. Let $D_+$ be any compact subset of $\Sigma_+\cap\mathcal{T}$. There exists $\eps_4>0$ such that for all $\eps\leq \eps_4, \ r^-_\eps(D_+)\subset D_-.$

\vspace{0.2in}
\noindent
{\em Proof.}  We prove (i). Let ${\bf x}=(w_0,\eta_0,\gamma(\eta_0))\in D_-$. Pick $\delta>0$ such that $U_{2\delta}=B_{2\delta}({\bf Z}^*_+)\cap\Sigma\subset D_+$.
Let $\psi_+$ denote the flow corresponding to subsystem \eqref{retA}-\eqref{retB}. Note $\psi_+((w_0,\eta_0),t)=(w(t),\eta(t))\rightarrow  (w^*_{+,s},\eta^*_{+,s})$ as $t\rightarrow\infty$. Additionally noting $\gamma(\eta)$ is continuous,  there exists  $T=T({\bf x})>0$ such that for all $t\geq T, $
$$ \|(w(t),\eta(t),\gamma(\eta(t)))-{\bf Z}^*_+\|<\delta.$$
Pick $\eps({\bf x})>0$ such that $t({\bf x},\eps({\bf x}))>T$ \ (recalling $\phi_+({\bf x}, t({\bf x},\eps({\bf x})))=r^+_{\eps({\bf x})}({\bf x})).$ Then for all $\eps\leq \eps({\bf x}), \ r^+_\eps({\bf x})\in  B_\delta({\bf Z}^*_+)\cap\Sigma.$

Let $c_2\in (0,1), \ c_2<\delta/\mbox{diam}(D_-)$. Pick $\eps_2>0$ as in Proposition 5.3(ii). Let $\eps_3=\min\{\eps({\bf x}), \eps_2\}$. For $\eps\leq \eps_3$ and ${\bf y}\in D_-$,
$$\|r^+_\eps({\bf y})-r^+_{\eps}({\bf x})\|\leq c_2\|{\bf y}-{\bf x}\|\leq c_2\mbox{ diam}(D_-)<\delta,$$
implying 
$$\|r^+_\eps({\bf y})- {\bf Z}^*_+\|\leq \|r^+_\eps({\bf y})-r^+_{\eps}({\bf x})\|+ \|r^+_{\eps}({\bf x})-{\bf Z}^*_+\|<2\delta.$$
Hence, $r^+_\eps(D_-)\subset U_{2\delta}\subset D_+$.  The proof of (ii) is similar. \ $_\square$

\vspace{0.2in}
\noindent
{\bf Theorem 5.6} \ With parameters as in Table 1, there exists $\hat{\eps}>0$ so that for all $\eps\leq \hat{\eps}$, system \eqref{Fil} admits an attracting periodic orbit.

\vspace{0.1in}
\noindent
{\bf Remark.} \ This result holds  more generally for any $0<b_0<b<b_1$, and for any $\tau>0$ and $\rho>0$.

\vspace{0.2in}
\noindent
{\em Proof.} Let $D_+\subset\Sigma_+\cap\mathcal{T}$, with $D_+$ compact and ${\bf Z}^*_+\in \mbox{Int}(D_+)$.  Let $D_-\subset\Sigma_-\cap\mathcal{T}$, with $D_-$ compact and ${\bf Z}^*_-\in \mbox{Int}(D_-)$.  Given $c_1\in (0,1)$, pick $\eps_1>0$ such that for all $\eps\leq \eps_1$ and for all ${\bf x}_1, {\bf x}_2\in D_+,$
$$\|r^-_\eps({\bf x}_2)-r^-_\eps({\bf x}_1)\|\leq c_1\|{\bf x}_2-{\bf x}_1\|.$$
Given $c_2\in (0,1)$, pick $\eps_2>0$ such that for all $\eps\leq \eps_2$ and for all ${\bf y}_1, {\bf y}_2\in D_-,$
$$\|r^+_\eps({\bf y}_2)-r^+_\eps({\bf y}_1)\|\leq c_2\|{\bf y}_2-{\bf y}_1\|.$$
Pick $\eps_3>0$ such that for all $\eps\leq \eps_3, \ r^+_\eps(D_-)\subset D_+,$ and choose $\eps_4>0$ such that for all $\eps\leq \eps_4, \ r^-_\eps(D_+)\subset D_-.$

For $\eps\leq \hat{\eps}=\min\{\eps_1,\eps_2,\eps_3,\eps_4\}, $
$$r_\eps=r^+_\eps\circ r^-_\eps :D_+\rightarrow D_+$$
is a contraction map with contraction factor $c_1c_2\in (0,1)$. Hence $r_\eps$ has a unique fixed point ${\bf x}^*$ to which all $r_\eps$-orbits in $D_+$ converge. Flowing via $\phi_-$ from ${\bf x}^*$ to $r^-_\eps({\bf x})={\bf y}^*$, and via $\phi_+$ from ${\bf y}^*$ to $r^+_\eps({\bf y}^*)={\bf x}^*$, provides the desired periodic orbit. Every Filippov trajectory of system \eqref{Fil} passing through $D_+$ converges to this limit cycle.
 \ $_\square$

\vspace{0.2in}
In Figure 5 we present the periodic behavior of the model for three $\eps$-values.  When $\eps$ is two orders of magnitude smaller than $\rho$, $w$ and $\eta$ move so quickly relative to $\xi$ that they essentially switch from a point near $(w^*_{+,s},\eta^*_{+,s})$ to a  point near $(w^*_{-,s},\eta^*_{-,s})$, and vice versa (Figure 5, top row). This is the essence of the idea behind the proof of Theorem 5.6.

The middle panel in Figure 5 illustrates the periodic behavior when $\eps$ is one order of magnitude smaller than $\rho$. Interestingly, the limit cycle  persists when $\eps$ and $\rho$ have the same order of magnitude (Figure 5, bottom row).

\begin{figure}[h!]\hspace*{0.3in}
\includegraphics[width=5.7in,trim = 1.6in 4.0in .8in  1.8in, clip]{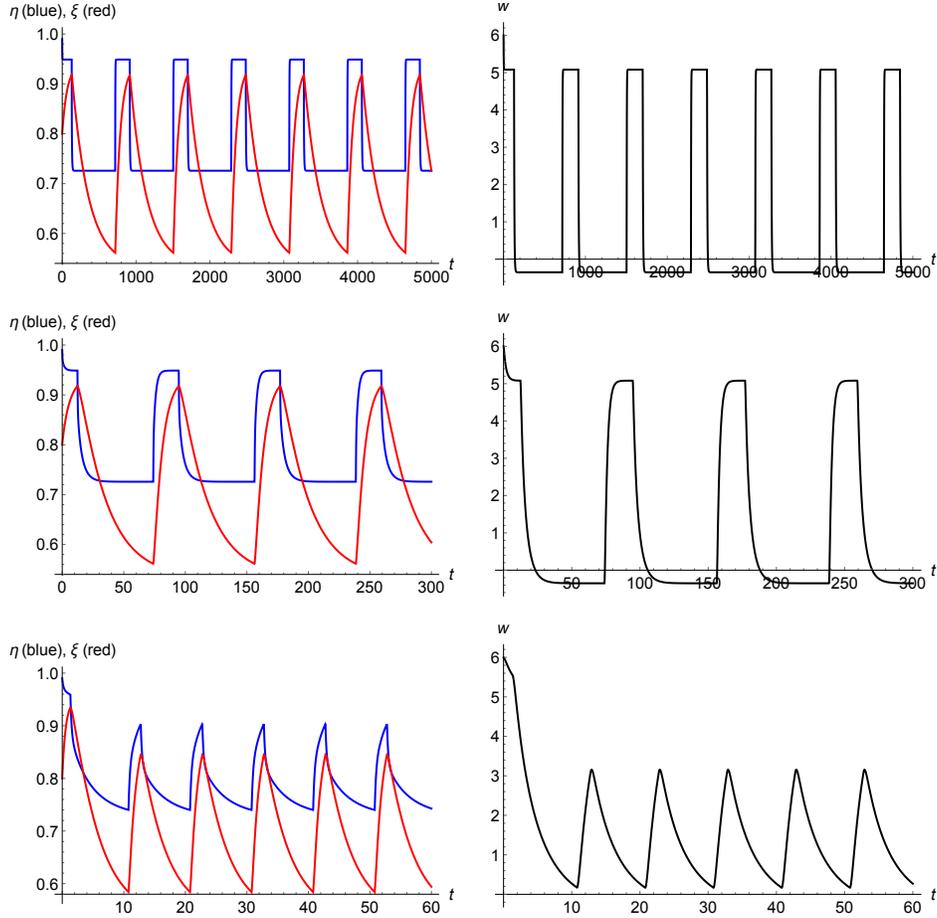} 
\begin{center}
\caption{{\small Periodic behavior for system \eqref{Fil}. {\em Top:} $\eps=0.003$.  {\em Middle:} $\eps=0.03$. \ \ {\em Bottom:} $\eps=0.3$.  In each plot $\tau=1$ and $ \rho=0.1$, with other parameters given in Table 1.}}
 \end{center}
\end{figure}

\section{Conclusion}

Building on previous work \cite{bud}, \cite{esther}, \cite{dickesther}, we separated the classic Budyko
ice line into two dynamic variables, the snow line and the glacial
boundary, producing a system of three ordinary differential equations.
Introducing two distinct climate states, one corresponding to glacial
advance and the other to glacial retreat, we formulated the equations
as a Filippov system with a two-dimensional plane forming the
discontinuity boundary.

For the parameter values of interest here, the system has two virtual
attracting equilibria, one for each of the two climate states.  For
example, starting with a climate where the glaciers are retreating,
the system moves toward its virtual equilibrium, crossing the
discontinuity boundary before reaching it.  Upon crossing the
boundary, the system switches climate regimes and starts heading
toward the virtual equilibrium for advancing glaciers.  Again, the
virtual equilibrium cannot be reached and the orbit crosses the
discontinuity boundary, switching back to retreating glaciers.  This
behavior is analogous to that found in the Welander model \cite{welander}, and it
produces a periodic orbit.

The proof of the existence of this periodic orbit uses ad-hoc singular
perturbation techniques, with $\epsilon$ as the singular parameter.
Recall that $\epsilon$ controls the speed at which the orbit
approaches the discontinuity boundary $\Sigma$.  By taking $\epsilon$
small enough, we were able to control the contraction parameter,
giving a unique attracting fixed point for the section map,
corresponding to an attracting periodic orbit for the system of
differential equations.  This periodic orbit represents an intrinsic
cycling of the climate systems, producing the glacial cycles.

\begin{figure}[h!]\hspace*{.35in}
\includegraphics[width=5.05in,trim = 1.4in 1.55in  1.2in 1.8in, clip]{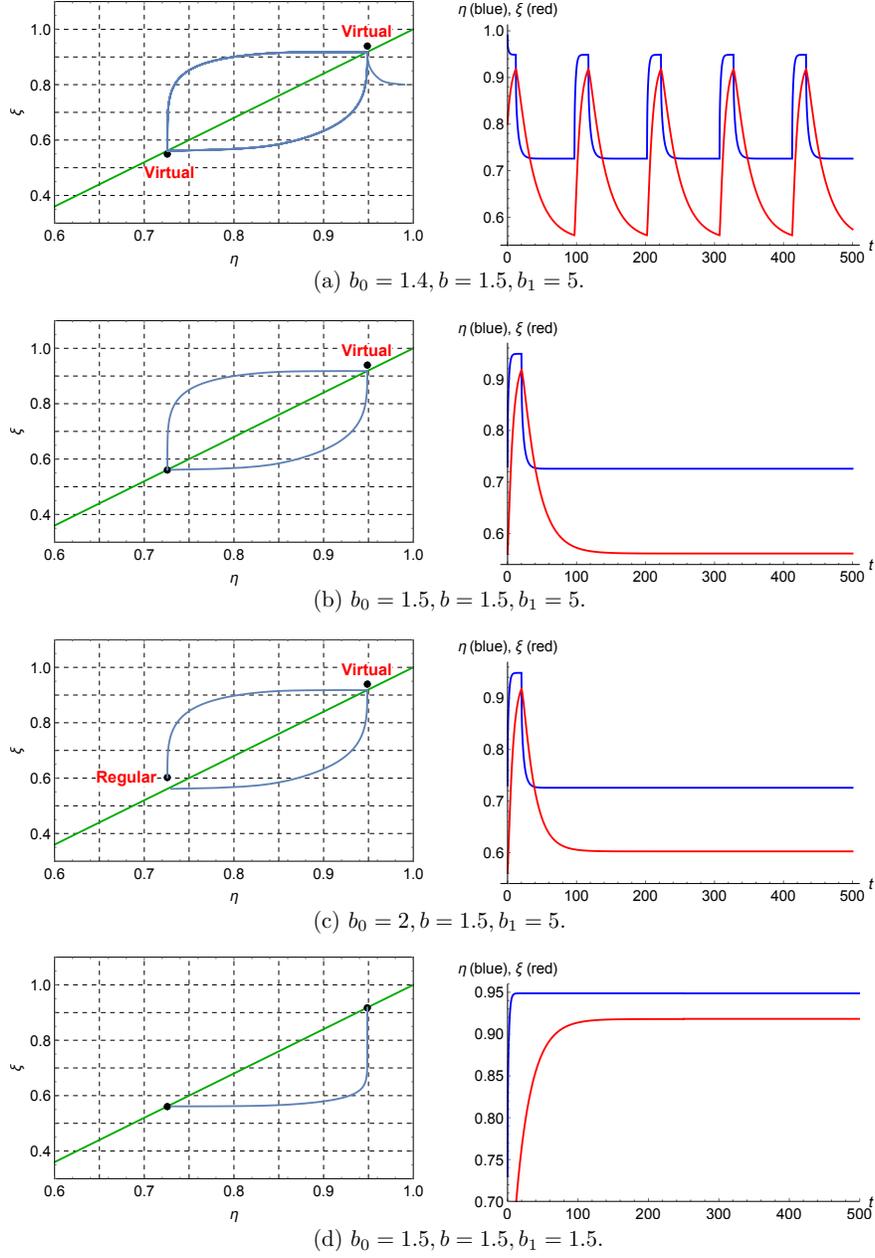} 

\vspace*{-.1in}
\begin{center}
\caption{{\small Illustration of a boundary equilibrium bifurcation as $b_0$ increases through $b$. Plotted are projections of various orbits onto the $(\eta,\xi)$-plane, and the corresponding $\eta(t)$ and $\xi(t)$ graphs. (a) Limit cycle. (b) ${\bf Q}^*_{-,s}\in\Sigma$.  (c) ${\bf Q}^*_{-,s}\in S_-$.  (d) ${\bf Q}^*_{-,s}, {\bf Q}^*_{+,s}\in \Sigma$.}}
\end{center}\end{figure}

We view our analysis as suggesting the need for a general theory of
singular perturbations for Filippov systems.

We also view our model as providing a first step toward a conceptual
model of the Earth's glacial cycles.  We chose parameters so that the
system has two virtual equilibria, giving us an intrinsic oscillation.
A next step would be to incorporate the Milankovitch cycles into the
model.  This would be easily accomplished, since the parameter $Q$ in
equation \eqref{budyko}  depends on the eccentricity while the parameter $s_2$ in
equation \eqref{approxsofy} depends on the obliquity \cite{dickclar}.  As these parameters vary
with the Milankovitch cycles, the system might be nudged across the
discontinuity boundary, creating a resonance with either eccentricity
or obliquity.

Other directions suggested by this model are more mathematical in
nature.  For example, both of the virtual equilibria are close to the
discontinuity boundary for the parameters we chose.  By adjusting the
parameters, one could produce a boundary equilibrium bifurcation for
one or both of the equilibria (see Figure 6). As one of the virtual equilibria
crossed the discontinuity boundary, the periodic orbit would
presumably disappear and be replaced by an attracting rest point.
Moving the parameters in the other direction would find the system
dramatically changing from an attracting equilibrium to a large
oscillation, reminiscent of the phenomenon used by Maasch and Salzman
to explain the mid-Pleistocene transition \cite{maasch}.

The mid-Pleistocene transition is a period roughly 1 million years ago
when the length of the glacial-interglacial cycles changed from
approximately 40,000 years to 100,000 years.  Substantial efforts have
been made to solve this ``100,000 year problem" (see, for example, \cite{ashwin}, \cite{huy2007}, \cite{paill1998}, \cite{parr2003}, \cite{tzip2003}).  One might allow the parameters to vary in
our model to see whether this transition can be realized.  Comparisons
can be made between model output and the paleoclimate data.

Again returning to mathematical questions, it would be interesting to
explore the parameter space with the goal of discovering the origin of
the periodic orbit we exhibited.  Perhaps the orbit originated with a
classic Hopf bifurcation, or perhaps it originated with bifurcation
from a ``fused focus" as discussed by Filippov \cite{fil}.  These explorations
may or may not have implications for interpretation of the
paleoclimate data, but they will almost certainly lead to interesting
examples of bifurcations in Filippov systems.

Conceptual in nature, the simple model introduced here gives rise to a
host of problems deserving of further study.

\vspace{0.3in}
\noindent
{\bf Acknowledgment.} The authors recognize and appreciate the support of the Mathematics and Climate Research Network (www.mathclimate.org). The work was partially supported by NSF Grants DMS-0940366 and
DMS-0940363.

\vspace{0.3in}
\renewcommand{\refname}

\end{document}